\documentclass[11pt]{article}

\usepackage[]{hyperref}

\usepackage{amsmath}
\usepackage{amssymb}
\usepackage{amsthm}
\usepackage{mathtools}
\usepackage{txfonts}

\newtheorem{theorem}{Theorem}
\newtheorem{lemma}[theorem]{Lemma}
\newtheorem{prop}[theorem]{Proposition}
\newtheorem{remark}{Remark}

\newcommand{\Ord}[1]{\mathcal{O}\left(#1\right)}

\newcommand{\eps}{\varepsilon}
\newcommand{\RR}{\mathbb{R}}
\newcommand{\NN}{\mathbb{N}}
\renewcommand{\rho}{\varrho}
\renewcommand{\phi}{\varphi}
\renewcommand{\theta}{\vartheta}
\newcommand{\mc}[1]{\mathcal{#1}}
\newcommand{\ut}{u_*}
\newcommand{\vt}{v_*}
\newcommand{\Olay}{\Omega_\mathrm{bl}}
\newcommand{\Oreg}{\Omega_\mathrm{reg}}

\newcommand{\E}{\mathrm{e}}
\newcommand{\D}{\mathrm{\, d}}
\newcommand{\abs}[1]{\left|#1\right|}
\newcommand{\norm}[1]{\left\|#1\right\|}
\newcommand{\enorm}[1]{\left|\left|\left|#1\right|\right|\right|}
\newcommand{\scal}[2]{\left\langle#1,#2\right\rangle}

\title{A balanced finite-element method for an\\ axisymmetrically loaded thin shell
\thanks{Supported by ANID-Chile through Fondecyt project 1230013}}

\author{Norbert Heuer\thanks{Facultad de Matemáticas, Pontificia Universidad
        Católica de Chile, Avenida Vicuña Mackenna 4860, Santiago, Chile
        (\texttt{nheuer@mat.uc.cl})}
        \and
        Torsten Lin\ss\thanks{{Fakult\"at f\"ur Mathematik und
        Informatik, FernUniversit\"at in Hagen, Universit\"atsstra{\ss}e 11, 58097 Hagen,
        Germany (\texttt{torsten.linss@fernuni-hagen.de})}}}
\date{}

\begin{document}
\maketitle

\begin{abstract}
  We analyse a finite-element discretisation of a differential equation
  describing an axisymmetrically loaded thin shell.
  The problem is singularly perturbed when the thickness of the shell becomes small.
  We prove robust convergence of the method in a balanced norm that captures the layers present
  in the solution. Numerical results confirm our findings.


\end{abstract}

\section{Introduction}

The deformation of thin elastic structures is a relevant subject in civil engineering,
and numerical methods are ubiquitous in the design of such structures.
The numerical analysis of thin-structure models is challenging as the presence of
singularities, layers, and locking phenomena is widespread, see, e.g.,
\cite{ChapelleB_11_FEA}.
In this paper we consider the axially symmetric deformation of a thin circular cylindrical shell,
as treated, e.g., in \cite[Section~5.5]{Fluegge_60_SS} and \cite[Chapter~15.3]{VentselK_01_TPS}.
This model is used, for example, in the design of pressurised cylindrical vessels like pipes and tanks.
The deformation is the solution of a singularly perturbed, fourth-order ordinary differential equation
and, depending on the boundary conditions, exhibits boundary layers. We consider the clamped
case, though other conditions like simply supported or free ends can be considered as well.
In dimensionless form the problem reads
\begin{subequations}\label{problem}
\begin{alignat}{2}
  \label{problem:ode}
  \eps^4 u^{(4)} + 4u & = f & \quad & \text{in} \ \ \Omega\coloneqq (0,1),\\
  \label{problem:bcs}
  u(\xi)=u'(\xi) & =0, & & \xi\in\{0,1\}.
\end{alignat}
\end{subequations}
Here, $f$ represents the load and \mbox{$0<\eps\ll 1$} is a small perturbation parameter.
It is proportional to the square root of the product of thickness times radius of the shell.
The equation also characterises the so-called simple edge effect
\cite{Goldenveizer_61_TET}, \cite[Chapter~17.5]{VentselK_01_TPS},
and appears in the Girkmann problem \cite{Girkmann_56_FEE}, see also
\cite{PitkaerantaBS_12_DRV,NiemiBPD_12_FEA,DevlooFGG_13_ACC,Niemi_16_BCS}.
For an early analysis of finite elements for thin circular shells we refer to
\cite{Morley_76_ADS} and the references cited there.

For constant coefficients, as specified here,
the solution to our model problem can be analytically represented.
This is precisely one of the main advantages of such model reductions.
In practice, however, computer simulations are standard and adapting analytical solutions to
boundary conditions is tedious.
Furthermore, numerical methods can be applied to structures of non-constant stiffness,
as in \cite{FrancisPT_90_BAA}. Analytical approaches have very limited flexibility,
cf. the early works by Olsson and Reissner \cite{Reissner_37_RTB,OlssonR_40_PBE}.
In conclusion, the use of numerical approaches for the solution of problem \eqref{problem}
is justified.

Main challenges in the numerical analysis of singularly perturbed problems are a lack of
stability (depending on the choice of norms) and poor approximations due to the presence
of boundary layers (depending on boundary conditions).
In the present case, the problem is uniformly stable with respect to the energy norm which,
on the other hand, does not control the layers of the solution.
We follow the strategy of Lin and Stynes \cite{LinS_12_BFE} for a reaction-dominated diffusion problem.
They derive a formulation that yields robust quasi-optimal convergence of a finite element
scheme in a stronger \emph{balanced} norm,
see also \cite{HeuerK_17_RDM} for an ultra-weak setting and
a related discontinuous Petrov--Galerkin (DPG) method.
Standard procedure to deal with the approximation of functions with boundary layers is to use
specific meshes, for instance Shishkin meshes \cite{Shishkin_92_DAS,MillerORS_96_FNM,Linss_10_LAM}.
This is also what we propose here, though our analysis can be extended to cover more general meshes
introduced in~\cite{RoosL_99_SCU}. In that way, we obtain a finite element
method that converges quasi-optimally in the balanced norm and which is robust with respect
to the perturbation parameter.

An overview of the remainder of this paper is as follows.
We start with collecting important properties of the solution to \eqref{problem}.
In \S\ref{sec_weak} we develop a variational formulation of problem \eqref{problem},
introduce the balance norm, and prove the unique solvability of our variational formulation
(Proposition~\ref{prop_well}). Afterwards, we introduce the finite element scheme with
Shishkin meshes and derive our main result of robust quasi-optimal convergence in the
balanced norm (Theorem~\ref{thm_conv} in \S\ref{sec_conv}). In \S\ref{sec_super} we
derive error estimates that prove superconvergence in the $L_2(\Omega)$-norm
(Theorem~\ref{thm_super}).
Finally, in \S\ref{sec_num} we report on some numerical experiments that confirm
the robust quasi-uniform convergence of our scheme and its superconvergence.

Let us finish this introduction by introducing some notation.
For $D\subset\Omega$, we consider
standard Lebesgue and Sobolev spaces $L_2(D)$ and $H^k(D)$ of integer order $k$,
respectively, and denote by $H^2_0(\Omega)\subset H^2(\Omega)$ the Sobolev space with
homogeneous boundary conditions. The $L_2(D)$-duality and norm are
$\scal{\cdot}{\cdot}_D$ and $\|\cdot\|_D$, respectively, and we drop
the index when $D=\Omega$. Furthermore, $C$ denotes a generic constant that is independent
of $\eps$ and the number of degrees of freedom of a discretisation.
Finally, notation ``$a\lesssim b$'' means $a\leq Cb$ with $C$ being a generic
positive constant, independent of any discretisation and perturbation parameter $\eps$.

\section{Properties of the exact solution} \label{sec_prop}

As is standard in finite element analysis, specific approximation properties for
problems with boundary layers require some knowledge of their behaviour.
This is what we provide in this section.

The layers in our problem are determined by the solutions of the
homogeneous problem
\begin{gather*}
   \eps^4 u^{(4)} + 4u = 0\,,
\end{gather*}
a fundamental system of which is given by
\begin{gather*}
  \exp\left(\pm \frac{x}{\eps}\right) \cos \frac{x}{\eps}
     \quad\text{and} \quad
  \exp\left(\pm \frac{x}{\eps}\right) \sin \frac{x}{\eps}\,.
\end{gather*}
The asymptotics of this problem were studied by O'Malley, see~\cite[\S3.B]{OMalley_91_SPA}.
Provided the right-hand side~$f$ of \eqref{problem} is sufficiently
smooth, the solution $u$ of~\eqref{problem} can be decomposed as
\begin{subequations}\label{decomp}
\begin{gather}
  \label{decomp:sum}
   u = u_{\mathrm{bl},0} + u_\mathrm{reg} + u_{\mathrm{bl},1}
          \quad \text{in} \ \bar\Omega,
  \intertext{where, for any given integer $K>0$, the regular component $u_\mathrm{reg}$ satisfies}
  \label{decomp:reg}
  \abs{u_\mathrm{reg}^{(\kappa)}(x)} \lesssim 1, \ \ x\in\bar\Omega,
     \ \ \kappa=0,1,\dots,K,
  \intertext{while for the layer components $u_{\mathrm{bl},0}$ and
  $u_{\mathrm{bl},1}$ the following bounds hold:}
  \label{decomp:lay}
  \abs{u_{\mathrm{bl},0}^{(\kappa)}(x)}
      \lesssim \eps^{-\kappa} \E^{- x/\eps}, \quad
  \abs{u_{\mathrm{bl},1}^{(\kappa)}(x)}
      \lesssim \eps^{-\kappa} \E^{-(1-x)/\eps},
                          \ \ x\in\bar\Omega, \ \ \kappa=0,1,\dots,K. \\
  \intertext{Moreover}
  \label{decomp:ode}
    \eps^4 u_\mathrm{reg}^{(4)} + 4u_\mathrm{reg} = f, \ \ \
    \eps^4 u_\mathrm{bl,i}^{(4)} + 4u_\mathrm{bl,i} = 0 \ \ \ \text{in} \ \ \Omega
      \ \ \ (i=0,1).
\end{gather}
\end{subequations}

\section{Variational formulation} \label{sec_weak}

The standard variational formulation of~\eqref{problem} reads:
Find $u\in H_0^2(\Omega)$ such that
\begin{gather*}
  a(u,v) \coloneqq
    \eps^4 \scal{u''}{v''} + 4 \scal{u}{v} = \scal{f}{v}\, \quad \forall v\in H_0^2(\Omega)\,.
\end{gather*}
The associated (standard) energy norm on $H^2(\Omega)$ is given by the square root of
\begin{gather*}
  \enorm{v}_s^2 \coloneqq a(v,v) = \eps^4 \norm{v''}^2 + 4 \norm{v}^2 \ \
        \forall v\in H^2(\Omega)\,.
\end{gather*}
A straightforward calculation shows that
\begin{gather*}
   \enorm{u_{\mathrm{bl},0}}_s + \enorm{u_{\mathrm{bl},1}}_s = \Ord{\eps^{1/2}}\,.
\end{gather*}
This means that the standard energy norm does not ``see'' the layers present
in the solution of~\eqref{problem}.
The correct weight for the norm of $u''$ would be $\eps^{3/2}$ instead
of $\eps^2$, in which case the resulting norm is called \emph{balanced}.
However, with respect to that norm the bilinear form $a(\cdot,\cdot)$ is not
uniformly coercive. The coercivity constant is of order $\eps$.

As indicated in the introduction, we adapt the technique from Lin and Stynes~\cite{LinS_12_BFE}
and define a variational formulation whose bilinear form induces a norm that is balanced.
The main idea is to test~\eqref{problem} with
$\eps^\alpha v^{(4)} + v$, $\alpha=3$, instead of $v$.
Note that a classical least-squares Galerkin method would use $\alpha=4$.

To avoid excessive $H^4$-regularity of both ansatz and test functions,
we reformulate~\eqref{problem} as a system of two second-order equations:
Find $u$ and $v$ such that
\begin{gather}\label{system}
   v -\eps^{3/2} u'' = 0 \quad \text{and} \quad
   \eps^{5/2} v'' + 4u  = f \quad \text{in} \ \Omega.
\end{gather}
Letting \mbox{$V \coloneqq H_0^2(\Omega)\times H^2(\Omega)$},
we choose a test function \mbox{$(\ut,\vt)\in V$} and
multiply the first equation with \mbox{$\vt -\eps^{3/2} \ut''$} and
the second one with \mbox{$\eps^{3/2} \vt'' + 4\ut$}.
We see that the solution $u$ of~\eqref{problem} satisfies
\begin{gather*}
   \scal{v -\eps^{3/2} u''}{\vt -\eps^{3/2} \ut''} = 0 \quad \text{and} \quad
   \scal{\eps^{5/2} v'' + 4u}{\eps^{3/2} \vt'' + 4\ut}
      = \scal{f}{\eps^{3/2} \vt'' + 4\ut}\,.
\end{gather*}
Next, let \mbox{$\lambda > 0$} be constant to be fixed later.
We define a bilinear form \mbox{$\mc{B}\colon V^2 \to \RR$} by
\begin{gather*}
  \mc{B}\bigl((u,v),(\ut,\vt)\bigr)
     \coloneqq \lambda \scal{v -\eps^{3/2} u''}{\vt -\eps^{3/2} \ut''}
               + \scal{\eps^{5/2} v'' + 4u}{\eps^{3/2} \vt'' + 4\ut}
\end{gather*}
and a linear functional $\mc{F} \colon V \to \RR$ by
\begin{gather*}
  \mc{F}\bigl((\ut,\vt)\bigr)
     \coloneqq \scal{f}{\eps^{3/2} \vt'' + 4\ut}.
\end{gather*}
Then a weak formulation of~\eqref{system} is: Find $(u,v)\in V$ such
that
\begin{gather}\label{weak}
  \mc{B}\bigl((u,v),(\ut,\vt)\bigr) = \mc{F}\bigl((\ut,\vt)\bigr) \quad
    \forall \ (\ut,\vt)\in V.
\end{gather}
We furnish $V$ we the norm (squared)
\begin{gather*}
  \enorm{(u,v)}^2 \coloneqq \norm{u}^2 + \eps^3\norm{u''}^2
    + \norm{v}^2 + \eps^4\norm{v''}^2\,.
\end{gather*}
Because of the $\eps^{3/2}\norm{u}$-term this norm is balanced.
Also note that the term $\eps^2\norm{v''}$ provides additional control
of the fourth-order derivative of $u$.

\emph{Notation:} Similarly as before, we shall use notation
$\mc{B}_D(\cdot,\cdot)$ and $\enorm{(\cdot)}_D$ when the integration in the
definitions above is restricted to $D\subset\Omega$.


\begin{lemma}\label{lem:coer-cont}
  If $\lambda\ge 3$ then bilinear form $\mc{B}$ is coercive and
  continuous with respect to the balanced norm $\enorm{\cdot}$,
  uniformly in $\eps$.
\end{lemma}
\begin{proof}
  First we prove coercivity of $\mc{B}$.
  Let $(u,v)\in V$ be arbitrary.
  Then
  \begin{align*}
    \mc{B}\bigl((u,v),(u,v)\bigr)
       & = \lambda \scal{v -\eps^{3/2} u''}{v -\eps^{3/2} u''}
                 + \scal{\eps^{5/2} v'' + 4u}{\eps^{3/2} v'' + 4u} \\
       & = \lambda \norm{v -\eps^{3/2} u''}^2
                 + \eps^4 \norm{v''}^2 + 16 \norm{u}^2 
                 + 4 \eps^{3/2}\left(1+\eps\right) \scal{v''}{u}.
  \end{align*}
  A direct calculation using $u\in H_0^2(\Omega)$ yields
  \begin{gather*}
    4 \eps^{3/2} \scal{v''}{u} =
    4 \eps^{3/2} \scal{v}{u''} = \norm{v+\eps^{3/2}u''}^2 - \norm{v-\eps^{3/2} u''}^2\,.
  \end{gather*}
  Thus
  \begin{align*}
    \mc{B}\bigl((u,v),(u,v)\bigr)
       & = \left(\lambda - 1 -\eps\right) \norm{v -\eps^{3/2} u''}^2
           + (1 +\eps) \norm{v + \eps^{3/2} u''}^2
                 + \eps^4 \norm{v''}^2 + 16 \norm{u}^2  \\
       & \ge \left(\lambda - 2 \right) \norm{v -\eps^{3/2} u''}^2
           + \norm{v + \eps^{3/2} u''}^2
                 + \eps^4 \norm{v''}^2 + 16 \norm{u}^2 ,
  \end{align*}
  where we have used that $\eps\in (0,1]$.
  Next, we note that
  \begin{gather*}
    \norm{v + \eps^{3/2} u''}^2 + \norm{v - \eps^{3/2} u''}^2
      = 2 \norm{v}^2 + 2 \eps^3 \norm{u''}^2\,.
  \end{gather*}
  Therefore,
  \begin{align*}
    \mc{B}\bigl((u,v),(u,v)\bigr)
       & \ge \left(\lambda - 3 \right) \norm{v -\eps^{3/2} u''}^2
           + 2 \norm{v}^2 + 2 \eps^3 \norm{u''}^2
                 + \eps^4 \norm{v''}^2 + 16 \norm{u}^2 \\
       & \ge 2 \norm{v}^2 + 2 \eps^3 \norm{u''}^2
                 + \eps^4 \norm{v''}^2 + 16 \norm{u}^2 \ge \enorm{(u,v)}^2\,.
  \end{align*}
  This is the coercivity of $\mc{B}(\cdot,\cdot)$.

  Next, we show the boundedness of $\mc{B}(\cdot,\cdot)$.
  Let $(u,v),(\ut,\vt)\in V$ be arbitrary.
  Integrating by parts, we get $\scal{u}{\vt''} = \scal{u''}{\vt}$.
  Hence,
  \begin{gather*}
    \begin{split}
      \mc{B}\bigl((u,v),(\ut,\vt)\bigr)
         & = \lambda \scal{v -\eps^{3/2} u''}{\vt -\eps^{3/2} \ut''} \\
         & \qquad\qquad
                   + \eps^4\scal{v''}{\vt''} + 4\eps^{5/2}\scal{v''}{\ut}
                   + 4\eps^{3/2}\scal{u''}{\vt} + 16\scal{u}{\ut}\,.
    \end{split}
  \end{gather*}
  Applying the Cauchy-Schwarz inequality and recalling the definition of
  $\enorm{\cdot}$ completes the proof.
\end{proof}

An immediate consequence of Lemma~\ref{lem:coer-cont},
the (non-uniform) boundedness of $\mc{F}$, and the Lax--Milgram
lemma is the following existence and uniqueness result.
\begin{prop} \label{prop_well}
  Let $f\in L_2(\Omega)$.
  If $\lambda\ge 3$ then~\eqref{weak} possesses a unique solution $(u,v)\in V$.
\end{prop}

\section{FEM discretisation on a layer-adapted mesh} \label{sec_fem}

Let $V_h=V_h^u\times V_h^v$ be a finite-dimensional subspace of $V$.
Then our discretisation reads: Find $(u_h,v_h)\in V_h$ such
that
\begin{gather}\label{fem}
  \mc{B}\bigl((u_h,v_h),(\ut,\vt)\bigr) = \mc{F}\bigl((\ut,\vt)\bigr) \quad
    \forall \ (\ut,\vt)\in V_h.
\end{gather}
Again Lemma~\ref{lem:coer-cont} and the Lax--Milgram lemma guarantee the
existence of a unique solution to~\eqref{fem}.
Furthermore, a standard argument proves that the method is quasi optimal.
\begin{prop}\label{prop:quasi-opt}
  Let $u$ be the solution of~\eqref{problem}. The approximation $(u_h,v_h)\in V_h$
  defined by~\eqref{fem} satisfies the robust, quasi-optimal error estimate
  \begin{align*}
    & \enorm{(u-u_h,\eps^{3/2} u''-v_h)}
      \le C \inf_{(\ut,\vt)\in V_h} \enorm{(u-\ut,\eps^{3/2} u''-\vt)}\\
    & \qquad \le C \inf_{\ut\in V_h^u}\left\{ \norm{u-\ut} + \eps^{3/2}\norm{(u-\ut)''}\right\}
          + C \eps^{3/2} \inf_{\vt\in V_h^v}\left\{ \norm{u''-\vt} + \eps^2\norm{u^{(4)} -\vt''}
                       \right\}.
  \end{align*}
\end{prop}

More specifically, we discretise~\eqref{weak} using conforming finite
elements with piecewise polynomials of degree $p\ge3$.
In the following, $\mc{P}_p$ denotes the space of polynomials of degree $p$.
For given $N\in \NN$ let
\begin{gather*}
  \omega\colon 0=x_0<x_1<\cdots<x_N=1
\end{gather*}
be an arbitrary partition of the domain \mbox{$\bar\Omega$}.
Set \mbox{$J_i\coloneqq\left(x_{i-1},x_i\right)$} and
\mbox{$h_i\coloneqq x_i - x_{i-1}$}, \mbox{$i=1,\dots,N$}.
For $p,k\in\NN_0$, we introduce the spline spaces
\begin{gather*}
  \mc{S}_p^k(\omega) \coloneqq \left\{ s\in C^k[0,1] \colon
    s|_{J_i} \in \mc{P}_p \ \forall i=1,\dots,N \right\}\,.
\end{gather*}
Then our finite element space is given by
\begin{gather*}
 V_h \coloneqq \mc{S}_p^1(\omega)^2 \cap V\,.
\end{gather*}
The boundary layers present in the solution of~\eqref{problem} may be
resolved using layer-adapted meshes.
As previously mentioned, we consider Shishkin meshes.

Let \mbox{$\sigma>0$} be a mesh parameter that will be fixed later and define
a mesh transition point by
\begin{gather}\label{def-tau}
  \tau \coloneqq \min \left\{\frac{1}{4}, \sigma\eps \ln N \right\}\,.
\end{gather}
Then the intervals $[0,\tau]$ and $[1-\tau,1]$ are uniformly dissected into
$N/4$ mesh intervals each, while $[\tau,1-\tau]$ is divided into $N/2$
subintervals.
In the sequel we shall restrict ourselves to the case {$\tau=\sigma\eps\ln N$}.
Otherwise the mesh is uniform and the argument proceeds in a standard manner
with $1/\eps \lesssim \ln N$ used in the derivative bounds of~\eqref{decomp}.

\subsection{A priori error estimate} \label{sec_conv}

In this section we state and prove the robust and quasi-optimal convergence of our
finite element scheme on Shishkin meshes. Before doing so, we need some preparation.
The convergence analysis starts from the quasi optimality, Proposition~\ref{prop:quasi-opt}.
It uses a specially designed function $(\tilde{u},\tilde{v})\in V_h$
that allows to bound the approximation error in terms of mesh parameter $N$.
The construction of $\tilde{u}$ and $\tilde{v}$ is based on the following
nodal interpolant. 

\begin{lemma}\label{lem:local-int}
  Let $p\in\NN_0$, $p>2k$, $J\coloneqq(a,b)$, $h\coloneqq b-a$ and $w\in H^{p+1}(I)$.
  \mbox{Let $I_Jw \in \mc{P}_p$} be defined by
  \begin{gather*}
    \left(I_J w- w\right)^{(\kappa)}(a)= 0\,, \ \ \kappa = 0,1,
    \intertext{and}
    \int_a^b \left(I_J w - w\right)''(x) q(x) \D x =0 \ \ \forall q\in\mc{P}_{p-2}\,.
  \end{gather*}
  Then
  \begin{gather*}
    \left(I_J w- w\right)^{(\kappa)}(b)= 0\,, \ \ \kappa = 0,1,
    \intertext{and}
      \norm{\left(I_J w- w\right)^{(k)}}_J
          \le C h^{p+1-\kappa} \norm{w^{(p+1)}}_J\,, \ \ \kappa = 0,1,2.
  \end{gather*}
\end{lemma}
\begin{proof}
  Note that \mbox{$\left(I_Jw\right)''\in\mc{P}_{p-2}$} is the $L_2$ projection
  of $w''$ onto $\mc{P}_{p-2}$.
  It can be defined equivalently as a truncated Legendre expansion of $w''$.
  Consequently, the technique from~\cite[\S3.3]{Schwab_98_php} can be applied to give
  the desired result.
\end{proof}

For a function \mbox{$w\in H^2(\Omega)$} we define a global
interpolant \mbox{$Iw\in\mc{S}_p^1(\omega)$} by applying the above definition
on each of the subintervals of our partition, i.e.,
\mbox{$I w|_{J_i} = I_{J_i}w$}, \mbox{$i=1,\dots,N$}.

In order to design our special representative
\mbox{$\left(\tilde{u},\tilde{v}\right)\in V_h$}
of \mbox{$\bigl(u,\eps^{3/2}u''\bigr)$},
we define auxiliary functions $\chi_0,\chi_1\in \mc{P}_3$ by
\begin{gather*}
  \chi_0(0)=\chi_0'(0)=\chi_0(\tau)=0, \quad \chi_0'(\tau)=1 \\
 \intertext{and}
  \chi_1(0)=\chi_1'(0)=\chi_1'(\tau)=0, \quad \chi_1(\tau)=1.
\end{gather*}
A direct calculation establishes the following bounds,
\begin{gather}\label{bound-chi}
  \norm{\chi_0}_{(0,\tau)} \lesssim \tau^{3/2}, \quad
  \norm{\chi_0''}_{(0,\tau)} \lesssim \tau^{-1/2}, \quad
  \norm{\chi_1}_{(0,\tau)} \lesssim \tau^{1/2}, \quad 
  \norm{\chi_1''}_{(0,\tau)} \lesssim \tau^{-3/2}\,.
\end{gather}
Recalling decomposition~\eqref{decomp}, let
\mbox{$u_{\mathrm{bl}} \coloneqq u_{\mathrm{bl},0} + u_{\mathrm{bl},1}$},
\mbox{$v_\mathrm{reg}\coloneqq \eps^{3/2} u_\mathrm{reg}''$} and
\mbox{$v_\mathrm{bl}\coloneqq \eps^{3/2} u_\mathrm{bl}''$}.
We define
\begin{align} \label{utilde}
  \tilde{u}(x)
    & \coloneqq
        \begin{cases}
          \left(I u\right)(x)
             - \chi_0(x) u_{\mathrm{bl}}'(\tau)
             - \chi_1(x) u_{\mathrm{bl}}(\tau)\,, & x\in [0,\tau], \\[.2ex]
          \bigl(I u_{\mathrm{reg}}\bigr)(x)\,, & x\in [\tau,1-\tau], \\[.2ex]
          \left(I u\right)(x)
             + \chi_0(1-x) u_{\mathrm{bl}}'(1-\tau)
             - \chi_1(1-x) u_{\mathrm{bl}}(1-\tau)\,, & x\in[1-\tau,1],
        \end{cases}
\end{align}
and
\begin{align} \label{vtilde}
  \tilde{v}(x)
    & \coloneqq
        \begin{cases}
          \left(I v\right)(x)
             - \chi_0(x) v_{\mathrm{bl}}'(\tau)
             - \chi_1(x) v_{\mathrm{bl}}(\tau)\,, & x\in [0,\tau], \\[.2ex]
          \bigl(I v_{\mathrm{reg}}\bigr)(x)\,, & x\in [\tau,1-\tau], \\[.2ex]
          \left(I v\right)(x)
             + \chi_0(1-x) v_{\mathrm{bl}}'(1-\tau)
             - \chi_1(1-x) v_{\mathrm{bl}}(1-\tau)\,, & x\in[1-\tau,1].
        \end{cases}
\end{align}
Note that
\begin{gather}\label{bl-in-tau}
  \abs{u_\mathrm{bl}(\xi)} \lesssim N^{-\sigma}, \ \ \
  \abs{u_\mathrm{bl}'(\xi)} \lesssim \eps^{-1} N^{-\sigma}, \ \ \
  \abs{v_\mathrm{bl}(\xi)} \lesssim \eps^{-1/2} N^{-\sigma}, \ \ \
  \abs{v_\mathrm{bl}'(\xi)} \lesssim \eps^{-3/2} N^{-\sigma}, \ \ \
   \xi\in\{\tau,1-\tau\}.
\end{gather}

\begin{lemma}\label{lem:u-ut}
  Let \mbox{$\Olay \coloneqq (0,\tau)\cup(1-\tau,1)$}
  and \mbox{$\Oreg \coloneqq (\tau,1-\tau)$}.
  Assume that~\eqref{decomp} holds for \mbox{$K=p+1$}.
  Then the bounds
  \begin{align} \label{u-ut:bl}
    \norm{u - \tilde{u}}_{\Olay}
       & \lesssim \eps^{1/2} \left\{\left(N^{-1} \ln N\right)^{p+1}
                         + N^{-\sigma} \left(\ln N\right)^{3/2}\right\} \,, \\
    \label{u-ut:bl''}
    \eps^{3/2} \norm{\left(u - \tilde{u}\right)''}_{\Olay}
       & \lesssim \left\{\left(N^{-1} \ln N\right)^{p-1}
          + N^{-\sigma}\right\} \,, \\
    \label{u-ut:reg}
    \norm{u - \tilde{u}}_{\Oreg}
       & \lesssim \eps^{1/2} N^{-\sigma} + N^{-(p+1)} \\
    \intertext{and}
    \label{u-ut:reg''}
    \eps^{3/2} \norm{u'' - \tilde{u}''}_{\Oreg}
       & \lesssim N^{-\sigma} + \eps^{3/2} N^{-(p-1)}
  \end{align}
  hold true.
\end{lemma}
\begin{proof}
  We study the two regions $\Olay$ and $\Oreg$ separately.
  Before, note that
  \begin{gather}\label{Hp}
    \norm{u^{(\kappa)}} \lesssim \eps^{-\kappa+1/2}
                            , \ \ \kappa=0,1,\dots,p+1, \ \ \text{by} \ \eqref{decomp}\,.
  \end{gather}

  \paragraph{$\Olay$:}
  We present the argument for $(0,\tau)$. Identical bounds hold for $(1-\tau,1)$.
  We have
  \begin{gather}\label{int-u-bl}
    \norm{u - \tilde{u}}_{(0,\tau)}
      \le \norm{u - I u}_{(0,\tau)}
        + \norm{\chi_0}_{(0,\tau)} \cdot \abs{u_{\mathrm{bl}}'(\tau)}
        + \norm{\chi_1}_{(0,\tau)} \cdot \abs{u_{\mathrm{bl}}(\tau)}\,.
  \end{gather}
  Lemma~\ref{lem:local-int} yields
  \begin{gather}\label{u-Iu}
    \norm{u - I u}_{(0,\tau)}^2
       = \sum_{i=1}^{N/4} \norm{u - I_{J_i} u}_{J_i}^2
       \lesssim \sum_{i=1}^{N/4} \left(\frac{\eps \ln N}{N}\right)^{2(p+1)}
                      \norm{u^{(p+1)}}_{J_i}^2 
       \le \left(\frac{\eps \ln N}{N}\right)^{2(p+1)}
                      \norm{u^{(p+1)}}^2\,.
  \end{gather}
  Thus
  \begin{gather*}
    \norm{u - I u}_{(0,\tau)}
       \lesssim \eps^{1/2} \left(N^{-1} \ln N\right)^{p+1}\,, \ \ \text{by} \ \eqref{Hp}.
  \end{gather*}
  Using~\eqref{bound-chi} and~\eqref{bl-in-tau}, we bound the remaining terms
  in~\eqref{int-u-bl} as follows:
  \begin{align*}
    \norm{\chi_0}_{(0,\tau)} \cdot \abs{u_{\mathrm{bl}}'(\tau)}
       + \norm{\chi_1}_{(0,\tau)} \cdot \abs{u_{\mathrm{bl}}(\tau)}
         \lesssim \eps^{1/2} \left(\ln N\right)^{3/2} N^{-\sigma}\,.
  \end{align*}
  This proves~\eqref{u-ut:bl}.

  Next, consider $(u-\tilde{u})''$,
  \begin{gather*}
    \norm{\left(u - \tilde{u}\right)''}_{(0,\tau)}
      \le \norm{\left(u - I u\right)''}_{(0,\tau)}
        + \norm{\chi_0''}_{(0,\tau)} \cdot \abs{u_{\mathrm{bl}}'(\tau)}
        + \norm{\chi_1''}_{(0,\tau)} \cdot \abs{u_{\mathrm{bl}}(\tau)}\,.
  \end{gather*}
  The arguments proceed in a similar manner, and we arrive at
  \begin{gather*}
    \norm{\left(u - I u\right)''}_{(0,\tau)}
         \lesssim \eps^{-3/2} \left(N^{-1} \ln N\right)^{p-1} \\
    \intertext{and}
    \norm{\chi_0''}_{(0,\tau)} \cdot \abs{u_{\mathrm{bl}}'(\tau)}
       + \norm{\chi_1''}_{(0,\tau)} \cdot \abs{u_{\mathrm{bl}}(\tau)}
         \lesssim \eps^{-3/2} N^{-\sigma}\,,
  \end{gather*}
  which gives~\eqref{u-ut:bl''}.

  \paragraph{$\Oreg$:}
  We have $u=u_\mathrm{bl} + u_\mathrm{reg}$.
  Therefore,
  \begin{gather*}
    \norm{u - \tilde{u}}_{\Oreg} \le
      \norm{u_\mathrm{bl}}_{\Oreg}
        + \norm{u_\mathrm{reg} - I u_\mathrm{reg}}_{\Oreg}\\
    \intertext{and}
    \eps^{3/2}\norm{\left(u - \tilde{u}\right)''}_{\Oreg} \le
      \eps^{3/2}\norm{u_\mathrm{bl}''}_{\Oreg}
        + \eps^{3/2}\norm{\bigl(u_\mathrm{reg} - I u_\mathrm{reg}\bigr)''}_{\Oreg}\,.
  \end{gather*}
  For $u_\mathrm{bl}$ we have
  \begin{gather*}
    \norm{u_\mathrm{bl}}_{\Oreg} \lesssim \eps^{1/2} N^{-\sigma}
       \ \ \ \text{and} \ \ \ 
    \norm{u_\mathrm{bl}''}_{\Oreg} \lesssim \eps^{-3/2} N^{-\sigma},
      \ \ \ \text{by~\eqref{decomp:lay} and~\eqref{def-tau}.}
  \end{gather*}
  The maximum mesh size of the Shishkin mesh is bounded by \mbox{$2/N$}.
  Therefore, Lemma~\ref{lem:local-int} and~\eqref{decomp:reg} yield
  \begin{gather*}
     \norm{u_{\mathrm{reg}} - I u_{\mathrm{reg}}}_{\Oreg}
           \lesssim N^{-(p+1)}
       \ \ \ \text{and} \ \ \ 
     \norm{\bigl(u_\mathrm{reg} - I u_\mathrm{reg}\bigr)''}_{\Oreg}
           \lesssim N^{-(p-1)}\,.
  \end{gather*}
  Combining the last four inequalities, we obtain~\eqref{u-ut:reg}
  and~\eqref{u-ut:reg''}.
\end{proof}

\begin{lemma}\label{lem:v-vt}
  Assume that~\eqref{decomp} holds for \mbox{$K=p+3$}.
  Then the bounds
  \begin{align} \label{v-vt:bl}
    \norm{v - \tilde{v}}_{\Olay}
       & \lesssim \left(N^{-1} \ln N\right)^{p+1}
                         + N^{-\sigma} \left(\ln N\right)^{3/2} \,, \\
    \label{v-vt:bl''}
    \eps^2 \norm{\left(v - \tilde{v}\right)''}_{\Olay}
       & \lesssim \left(N^{-1} \ln N\right)^{p-1} + N^{-\sigma} \,, \\
    \label{v-vt:reg}
    \norm{v - \tilde{v}}_{\Oreg}
       & \lesssim N^{-\sigma} + N^{-(p+1)} \\
    \intertext{and}
    \label{v-vt:reg''}
    \eps^2 \norm{v'' - \tilde{v}''}_{\Oreg}
       & \lesssim N^{-\sigma} + \eps^2 N^{-(p-1)}
  \end{align}
  hold true.
\end{lemma}
\begin{proof}
  The arguments are very similar to the ones used in the proof of Lemma~\ref{lem:u-ut}.

  \paragraph{$\Olay$:}
  We have
  \begin{gather}\label{int-v-bl}
    \norm{v - \tilde{v}}_{(0,\tau)}
      \le \norm{v - I v}_{(0,\tau)}
        + \norm{\chi_0}_{(0,\tau)} \cdot \abs{v_{\mathrm{bl}}'(\tau)}
        + \norm{\chi_1}_{(0,\tau)} \cdot \abs{v_{\mathrm{bl}}(\tau)}\,.
  \end{gather}
  By Lemma~\ref{lem:local-int} we obtain
  \begin{gather*}
    \norm{v - I v}_{(0,\tau)}^2
       = \sum_{i=1}^{N/4} \norm{v - I_{J_i} v}_{J_i}^2
       \lesssim \eps^3 \sum_{i=1}^{N/4} \left(\frac{\eps \ln N}{N}\right)^{2(p+1)}
                      \norm{u^{(p+3)}}_{J_i}^2
       \le \eps^3 \left(\frac{\eps \ln N}{N}\right)^{2(p+1)}
                      \norm{u^{(p+3)}}^2\,.
  \end{gather*}
  Thus
  \begin{gather}\label{v-Iv}
    \norm{v - I v}_{(0,\tau)} \lesssim \left(N^{-1} \ln N\right)^{p+1}\,,
     \quad \text{by~\eqref{Hp}.}
  \end{gather}
  Using~\eqref{bound-chi} and~\eqref{bl-in-tau}, the remaining terms
  in~\eqref{int-v-bl} are bounded as follows:
  \begin{align*}
   \norm{\chi_0}_{(0,\tau)} \cdot \abs{v_{\mathrm{bl}}'(\tau)}
     + \norm{\chi_1}_{(0,\tau)} \cdot \abs{u_{\mathrm{bl}}(\tau)}
     & \lesssim \left(\ln N\right)^{3/2} N^{-\sigma}\,.
  \end{align*}
  This yields~\eqref{v-vt:bl}.

  Next, we consider $(v-\tilde{v})''$,
  \begin{gather*}
    \norm{\left(v - \tilde{v}\right)''}_{(0,\tau)}
      \le \norm{\left(v - I v\right)''}_{(0,\tau)}
        + \norm{\chi_0''}_{(0,\tau)} \cdot \abs{v_{\mathrm{bl}}'(\tau)}
        + \norm{\chi_1''}_{(0,\tau)} \cdot \abs{v_{\mathrm{bl}}(\tau)}\,.
  \end{gather*}
  The argument proceeds as before, yielding
  \begin{gather*}
    \norm{\left(v - I v\right)''}_{(0,\tau)}
       \lesssim \eps^{-2} \left(N^{-1} \ln N\right)^{p-1} \\
    \intertext{and}
    \norm{\chi_0''}_{(0,\tau)} \cdot \abs{v_{\mathrm{bl}}'(\tau)}
     + \norm{\chi_1''}_{(0,\tau)} \cdot \abs{v_{\mathrm{bl}}(\tau)}
     \lesssim \eps^{-2} N^{-\sigma}\,.
  \end{gather*}
  We conclude that~\eqref{v-vt:bl''} holds true.

  \paragraph{$\Oreg$:}
  We have $v=v_\mathrm{bl} + v_\mathrm{reg}$.
  Therefore,
  \begin{gather*}
    \norm{v - \tilde{v}}_{\Oreg} \le
      \norm{v_\mathrm{bl}}_{\Oreg}
        + \norm{v_\mathrm{reg} - I v_\mathrm{reg}}_{\Oreg}\\
    \intertext{and}
    \eps^{3/2}\norm{\left(v - \tilde{v}\right)''}_{\Oreg} \le
      \eps^{3/2}\norm{v_\mathrm{bl}''}_{\Oreg}
        + \eps^{3/2}\norm{\bigl(v_\mathrm{reg} - I v_\mathrm{reg}\bigr)''}_{\Oreg}\,.
  \end{gather*}
  For $v_\mathrm{bl}$ we have
  \begin{gather*}
    \norm{v_\mathrm{bl}}_{\Oreg}
      = \eps^{3/2} \norm{u_\mathrm{bl}''}_{\Oreg}
      \lesssim N^{-\sigma}
       \ \ \ \text{and} \ \ \ 
    \norm{v_\mathrm{bl}''}_{\Oreg}
      = \eps^{3/2} \norm{u_\mathrm{bl}^{(4)}}_{\Oreg}
      \lesssim \eps^{-2} N^{-\sigma},
      \ \ \ \text{by~\eqref{decomp:lay} and~\eqref{def-tau}.}
  \end{gather*}
  The maximum mesh size is bounded by $2/N$.
  Therefore, Lemma~\ref{lem:local-int} and~\eqref{decomp:reg} yield
  \begin{gather*}
     \norm{v_{\mathrm{reg}} - I v_{\mathrm{reg}}}_{\Oreg}
           \lesssim N^{-(p+1)}
       \ \ \ \text{and} \ \ \ 
     \norm{\bigl(v_\mathrm{reg} - I v_\mathrm{reg}\bigr)''}_{\Oreg}
           \lesssim N^{-(p-1)}\,.
  \end{gather*}
  Combining the last four inequalities, we obtain~\eqref{v-vt:reg}
  and~\eqref{v-vt:reg''}.
\end{proof}

We are ready to prove our first main result.

\begin{theorem} \label{thm_conv}
  We use bilinear form $\mc{B}(\cdot,\cdot)$ with $\lambda\ge 3$
  and Shishkin meshes with parameters $p\ge3$ and $\sigma\ge p-1$.
  Assume that \eqref{decomp} holds for \mbox{$K=p+3$}.
  Then scheme~\eqref{fem} is uniformly convergent with
  \begin{gather*}
    \enorm{(u-u_h,\eps^{3/2} u''-v_h)}
       \lesssim \left(N^{-1} \ln N\right)^{p-1}\,.
  \end{gather*}
\end{theorem}
\begin{proof}
  By construction of \mbox{$\tilde{u} \in V_h^u$}, using Lemma~\ref{lem:u-ut} we can bound
  \begin{gather*}
    \begin{split}
    & \inf_{\ut\in V_h^u}\left\{ \norm{u-\ut} + \eps^{3/2}\norm{(u-\ut)''}\right\} \\
    & \qquad\qquad
      \le \norm{u - \tilde{u}}_{\Olay}
        + \eps^{3/2} \norm{\left(u - \tilde{u}\right)''}_{\Olay}
        + \norm{u - \tilde{u}}_{\Oreg}
        + \eps^{3/2} \norm{\left(u - \tilde{u}\right)''}_{\Oreg}\,. \\
    & \qquad\qquad
      \lesssim \left(N^{-1} \ln N\right)^{p-1},
    \end{split}
  \end{gather*}
  Similarly, employing Lemma~\ref{lem:v-vt}, we have
  \begin{gather*}
    \inf_{\vt\in V_h^v}\left\{\norm{v''-\vt} + \eps^2\norm{v'' -\vt''}\right\}
         \lesssim \left(N^{-1} \ln N\right)^{p-1}\,.
  \end{gather*}
  An application of Proposition~\ref{prop:quasi-opt} proves the statement of the
  theorem.
\end{proof}


\subsection{Superconvergence} \label{sec_super}

In this section we establish our second main result, the superconvergence of scheme \eqref{fem}
on Shishkin meshes. We start by proving a superconvergence property.

\begin{prop}\label{prop_super}
  Let $f\in L_2(\Omega)$ be given.
  We use bilinear form $\mc{B}(\cdot,\cdot)$ with $\lambda\ge 3$ and 
  Shishkin meshes with parameters $p\ge3$ and $\sigma\ge p+1$.
  Assume that~\eqref{decomp} holds for \mbox{$K=p+3$}.
  Then scheme~\eqref{fem} is uniformly superconvergent with
  \begin{gather*}
    \enorm{(\tilde{u}-u_h, \tilde{v}-v_h)}
      \lesssim \left\{\left(\ln N\right)^{p+1} + \abs{\lambda-4}
            \min\left\{\eps^{3/2}N^{2},\eps^{-1/2}\right\} \right\}
         N^{-(p+1)}\,.
  \end{gather*}
  Here, $\tilde{u}$ and $\tilde{v}$ are the functions defined in
  \eqref{utilde} and \eqref{vtilde}, respectively.
\end{prop}

\begin{proof}
Let \mbox{$\ut\coloneqq \tilde{u}-u_h$} and \mbox{$\vt\coloneqq \tilde{v}-v_h$}.
Coercivity of $\mc{B}(\cdot,\cdot)$ and Galerkin orthogonality imply
\begin{gather}\label{super-start}
\begin{split}
  \enorm{\left(\ut,\vt\right)}^2
     & \le \mc{B} \bigl(\left(\ut,\vt\right),\left(\ut,\vt\right)\bigr)
       = \mc{B} \bigl(\left(\tilde{u}-u,\tilde{v}-v\right),\left(\ut,\vt\right)\bigr) \\
     & = \lambda\scal{\tilde{v}-v -\eps^{3/2}\tilde{u}'' + \eps^{3/2}u''}{\vt-\eps^{3/2}\ut''}
          + \scal{\eps^{5/2}(\tilde{v}-v)'' + 4(\tilde{u}-u)}{\eps^{3/2}\vt''+4\ut}\,.
\end{split}
\end{gather}
Consider the layer region $(0,\tau)$.
The definition of interpolant $I$ and integration by parts show that
\begin{align*}
  & \scal{\tilde{v}-v -\eps^{3/2}\tilde{u}'' + \eps^{3/2}u''}{\vt-\eps^{3/2}\ut''}_{(0,\tau)} \\
  & \qquad =
    \scal{-\chi_0v_\mathrm{bl}'(\tau) - \chi_1 v_\mathrm{bl}(\tau)
          +\eps^{3/2}\chi_0''u_\mathrm{bl}'(\tau) + \eps^{3/2} \chi_1'' u_\mathrm{bl}(\tau)}{\vt-\eps^{3/2}\ut''}_{(0,\tau)} \\
  & \qquad \qquad
      + \scal{Iv-v}{\vt-\eps^{3/2}\ut''}_{(0,\tau)}
      - \eps^{3/2}\scal{Iu-u}{\vt''}_{(0,\tau)} \\
 \intertext{and}
  & \scal{\eps^{5/2}(\tilde{v}-v)'' + 4(\tilde{u}-u)}{\eps^{3/2}\vt''+4\ut}_{(0,\tau)} \\
  & \qquad =
    - \scal{\eps^{5/2}\left(\chi_0''v_\mathrm{bl}'(\tau) + \chi_1'' v_\mathrm{bl}(\tau)\right)
          + 4\left(\chi_0u_\mathrm{bl}'(\tau) + \chi_1 u_\mathrm{bl}(\tau)\right)}{\eps^{3/2}\vt''+4\ut}_{(0,\tau)} \\
  & \qquad \qquad
          + 4\eps^{3/2} \scal{Iu-u}{\vt''}_{(0,\tau)} + 4\eps^{5/2} \scal{Iv-v}{\ut''}_{(0,\tau)} + 16 \scal{Iu-u}{\ut}_{(0,\tau)}\,.
\end{align*}
Bounds~\eqref{bound-chi} and~\eqref{bl-in-tau} give
\begin{align*}
  & \abs{\scal{-\chi_0v_\mathrm{bl}'(\tau) - \chi_1 v_\mathrm{bl}(\tau)
          +\eps^{3/2}\chi_0''u_\mathrm{bl}'(\tau) + \eps^{3/2} \chi_1'' u_\mathrm{bl}(\tau)}{\vt-\eps^{3/2}\ut''}_{(0,\tau)}} \\
  & \qquad\quad \lesssim \left(\ln N\right)^{3/2} N^{-\sigma} \norm{\vt-\eps^{3/2}\ut''}_{(0,\tau)}
                \lesssim \left(\ln N\right)^{3/2} N^{-\sigma} \enorm{(\ut,\vt)}_{(0,\tau)} \\
  \intertext{and}
  & \abs{\scal{\eps^{5/2}\left(\chi_0''v_\mathrm{bl}'(\tau) + \chi_1'' v_\mathrm{bl}(\tau)\right)
          + 4\left(\chi_0u_\mathrm{bl}'(\tau) + \chi_1 u_\mathrm{bl}(\tau)\right)}{\eps^{3/2}\vt''+4\ut}_{(0,\tau)}} \\
  & \qquad\quad \lesssim \eps^{1/2} \left(\ln N\right)^{3/2} N^{-\sigma} \norm{\eps^{3/2}\vt''+4\ut}_{(0,\tau)}
                \lesssim \left(\ln N\right)^{3/2} N^{-\sigma} \enorm{(\ut,\vt)}_{(0,\tau)}\,.
\end{align*}
For the terms involving interpolation we have
\begin{align*}
   \eps^{3/2} \abs{\scal{Iu-u}{\vt''}_{(0,\tau)}}
     + \abs{\scal{Iu-u}{\ut}_{(0,\tau)}}
      & \lesssim \left(N^{-1} \ln N\right)^{p+1}\enorm{(\ut,\vt)}_{(0,\tau)}
               \,, \quad \text{by~\eqref{u-Iu},} \\
  \intertext{and}
   \abs{\scal{Iv-v}{\vt-\eps^{3/2}\ut''}_{(0,\tau)}}
   + \eps^{5/2} \abs{\scal{Iv-v}{\ut''}_{(0,\tau)}}
      & \lesssim \left(N^{-1} \ln N\right)^{p+1}\enorm{(\ut,\vt)}_{(0,\tau)}
               \,, \quad \text{by~\eqref{v-Iv}.}
\end{align*}
Gathering all these estimates and using that $\sigma\ge p+1$, we obtain
\begin{gather}\label{super-bl}
  \abs{\mc{B}_{(0,\tau)} \bigl(\left(\tilde{u}-u,\tilde{v}-v\right),\left(\ut,\vt\right)\bigr)}
    \lesssim \left(N^{-1} \ln N\right)^{p+1}\enorm{(\ut,\vt)}_{(0,\tau)}\,.
\end{gather}
Clearly, an identical bound holds for the layer region $(1-\tau,1)$.

On $(\tau,1-\tau)$, we have $\ut=Iu_\mathrm{reg}$ and $\vt=Iv_\mathrm{reg}$.
Integration by parts, the definition of interpolation operator $I$ and
\eqref{decomp:ode} give
\begin{gather}\label{super-aux}
\begin{split}
  & \lambda \scal{\tilde{v}-v -\eps^{3/2}\tilde{u}'' + \eps^{3/2}u''}{\vt-\eps^{3/2}\ut''}_{(\tau,1-\tau)}
  + \scal{\eps^{5/2}(\tilde{v}-v)'' + 4(\tilde{u}-u)}{\eps^{3/2}\vt''+4\ut}_{(\tau,1-\tau)} \\
  & \qquad =
    \lambda \scal{Iv_\mathrm{reg} - v_\mathrm{reg}}{\vt-\eps^{3/2}\ut''}_{(\tau,1-\tau)}
      + (4 - \lambda) \eps^{3/2}\scal{Iu_\mathrm{reg}-u_\mathrm{reg}}{\vt''}_{(\tau,1-\tau)} \\
  & \qquad\qquad\qquad +
           4\eps^{5/2} \scal{Iv_\mathrm{reg}-v_\mathrm{reg}}{\ut''}_{(\tau,1-\tau)}
         + 16 \scal{Iu_\mathrm{reg}-u_\mathrm{reg}}{\ut}_{(\tau,1-\tau)}
\end{split}
\end{gather}
Using Lemma~\ref{lem:local-int} and~\eqref{decomp:reg}, we obtain the following
$L_2$-norm bounds for the interpolation error:
\begin{gather*}
    \norm{Iv_\mathrm{reg} - v_\mathrm{reg}}_{(\tau,1-\tau)}
    + \norm{Iu_\mathrm{reg} - u_\mathrm{reg}}_{(\tau,1-\tau)} \lesssim N^{-(p+1)}\,.
\end{gather*}
An application of the Cauchy-Schwarz inequality gives
\begin{gather}\label{super-reg-1}
\begin{split}
  & \abs{\lambda \scal{Iv_\mathrm{reg} - v_\mathrm{reg}}{\vt-\eps^{3/2}\ut''}_{(\tau,1-\tau)}
         + 4\eps^{5/2} \scal{Iv_\mathrm{reg}-v_\mathrm{reg}}{\ut''}_{(\tau,1-\tau)}
         + 16 \scal{Iu_\mathrm{reg}-u_\mathrm{reg}}{\ut}_{(\tau,1-\tau)}} \\
  & \qquad
     \lesssim N^{-(p+1)} \enorm{(\ut,\vt)}_{(\tau,1-\tau)}
\end{split}
\end{gather}
and
\begin{gather*}
  \abs{(4 - \lambda) \eps^{3/2}\scal{Iu_\mathrm{reg}-u_\mathrm{reg}}{\vt''}_{(\tau,1-\tau)}}
      \lesssim
  \abs{4 - \lambda} N^{-(p+1)} \eps^{3/2} \norm{\vt''}_{(\tau,1-\tau)}\,.
\end{gather*}
Note that
\begin{align*}
  \eps^{3/2} \norm{\vt''}_{(\tau,1-\tau)} & \lesssim
  \eps^{-1/2} \enorm{(\ut,\vt)}_{(\tau,1-\tau)}\,,
  \intertext{by the definition of the balanced norm, and that}
  \eps^{3/2} \norm{\vt''}_{(\tau,1-\tau)} & \lesssim
  \eps^{3/2} N^2 \norm{\vt}_{(\tau,1-\tau)}
  \lesssim \eps^{3/2} N^2 \enorm{(\ut,\vt)}_{(\tau,1-\tau)}\,,
\end{align*}
by an inverse inequality.
Thus
\begin{align*}
  & \abs{(4 - \lambda) \eps^{3/2}\scal{Iu_\mathrm{reg}-u_\mathrm{reg}}{\vt''}_{(\tau,1-\tau)}}
      \lesssim
  \abs{4 - \lambda} N^{-(p+1)} \min\left\{\eps^{3/2} N^2,\eps^{-1/2}\right\}
              \enorm{(\ut,\vt)}_{(\tau,1-\tau)}\,,
\end{align*}
and
\begin{gather}\label{super-reg}
  \abs{\mc{B}_{(\tau,1-\tau)} \bigl(\left(\tilde{u}-u,\tilde{v}-v\right),\left(\ut,\vt\right)\bigr)}
    \lesssim \left\{1 + \abs{\lambda-4} \min\left\{\eps^{3/2} N^2,\eps^{-1/2}\right\}\right\}
         N^{-(p+1)}\enorm{(\ut,\vt)}_{(\tau,1-\tau)}\,.
\end{gather}
Finally, combining~\eqref{super-start}, \eqref{super-bl} and~\eqref{super-reg},
we obtain the statement of the theorem.
\end{proof}

\begin{remark}
  If one chooses $\lambda=4$ or if $\eps^{3/2}\le N^{-2}$, which is typical for singularly
  perturbed problems, then the error bound simplifies to
  \begin{gather*}
    \enorm{(\tilde{u}-u_h, \tilde{v}-v_h)}
      \lesssim \left(N^{-1} \ln N\right)^{(p+1)}\,.
  \end{gather*}
\end{remark}

Consequence of the superconvergence property established by Proposition~\ref{prop_super},
combined with an application of Lemmas~\ref{lem:u-ut},~\ref{lem:v-vt} and the triangle inequality,
is the superconvergence of scheme \eqref{fem} in the $L_2(\Omega)$-norm.

\begin{theorem} \label{thm_super}
  Let the assumptions of Proposition~\ref{prop_super} hold true.
  Then
  \begin{gather*}
    \norm{u-u_h} + \norm{\eps^{3/2} u'' - v_h}
      \lesssim \left\{\left(\ln N\right)^{p+1} + \abs{\lambda-4}
       \min\left\{\eps^{3/2} N^2,\eps^{-1/2}\right\}\right\}
         N^{-(p+1)}\,.
  \end{gather*}
  If $\lambda=4$ or $\eps^{3/2}\le N^{-2}$ then
  \begin{gather*}
    \norm{u-u_h} + \norm{\eps^{3/2} u'' - v_h}
      \lesssim \left(N^{-1} \ln N\right)^{(p+1)}\,.
  \end{gather*}
\end{theorem}

The superconvergence is indeed observed in our numerical experiments reported next.

\section{Numerical results} \label{sec_num}

We consider problem \eqref{problem} with given solution
$u(x)=\E^{-x/\varepsilon} \cos(x/\varepsilon) + \E^x$,
corresponding right-hand side function $f$ and essential boundary data $u(0),u'(0),u(1),u'(1)$.
Function $u$ has only a boundary layer at $x=0$ so that our numerical scheme uses
Shishkin meshes refined towards $0$, with transition point $\tau=\min\{1/2,\sigma\varepsilon\ln N\}$
($\sigma:=4$) and $N/2$ elements both in $[0,\tau]$ and $[\tau,1]$.
We use polynomial degree $p=3$ throughout, both for $u_h$ and $v_h$.

\begin{figure}
\begin{center}
\includegraphics[width=0.8\textwidth]{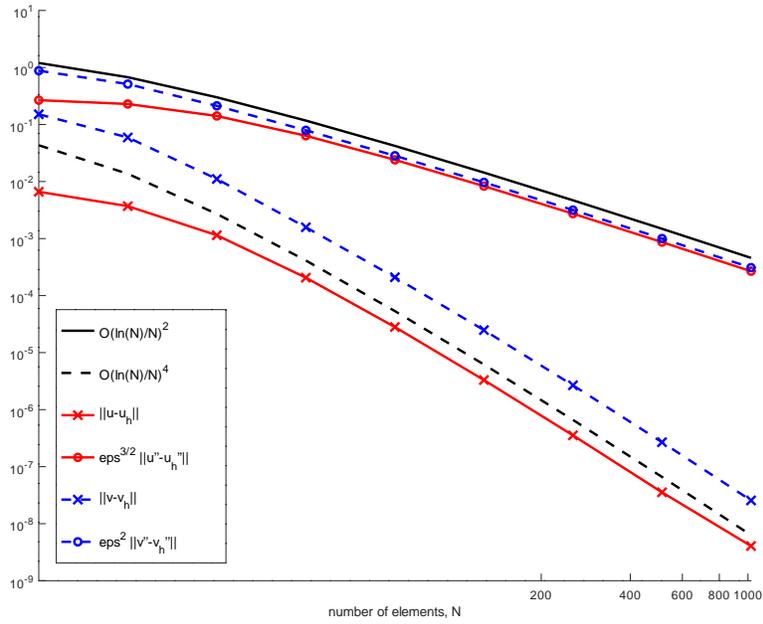}
\end{center}
\caption{Shishkin meshes with $N\in [4,1024]$ elements, $\mathrm{eps}=\varepsilon=10^{-2}$.}
\label{fig_N}
\end{figure}

Figure~\ref{fig_N} shows the individual terms of the error in balanced norm versus the number
of elements $N$. They are $\|u-u_h\|$, $\eps^{3/2}\|u''-u_h''\|$, $\|v-v_h\|$,
and $\eps^2\|v''-v_h''\|$ with $v=\eps^{3/2}u''$.
Also curves of orders $(N^{-1}\ln N)^2$ and $(N^{-1}\ln N)^4$ are shown, and note that
we have used logarithmic scales for both axes. We observe convergence orders
$\|u-u_h\|, \|v-v_h\|=O(N^{-1}\ln N)^4$, that is, superconvergence in $L_2(\Omega)$-norm, and
$\eps^{3/2}\|u''-u_h''\|, \eps^2\|v''-v_h''\|=O(N^{-1}\ln N)^2$.
These are the expected orders stated by Theorems~\ref{thm_conv} and~\ref{thm_super}.
We remark that we did not observe a lack of superconvergence for tests with $\lambda\not=4$.

\begin{figure}
\begin{center}
\includegraphics[width=0.8\textwidth]{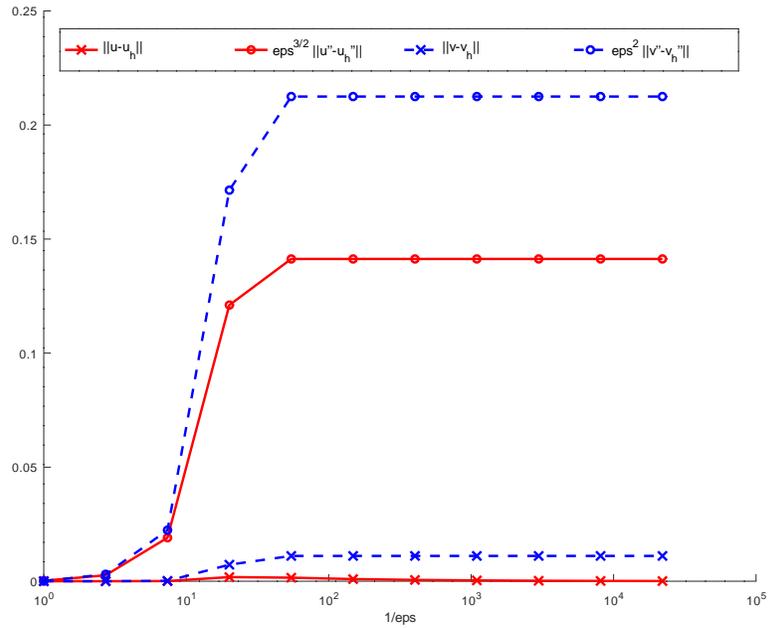}
\end{center}
\caption{Shishkin mesh with $N=16$ elements, $\mathrm{eps}=\varepsilon\in [1,\E^{-10}]$.}
\label{fig_eps}
\end{figure}

Robustness of our method is illustrated by Figure~\ref{fig_eps}. It shows
the same error terms from before, now versus $1/\eps$ and using a logarithmic scale only for the abscissa.
We consider Shishkin meshes with a fixed number of $N=16$ elements,
and vary $\eps$ between $1$ and $\E^{-10}$. All the individual error terms
quickly tend to small constants for decreasing $\eps$,
thus confirming the a priori error estimate by Theorem~\ref{thm_conv}
with hidden constant independent of $\eps$.

\begin{thebibliography}{10}

\bibitem{ChapelleB_11_FEA}
{\sc D.~Chapelle and K.-J. Bathe}, {\em The finite element analysis of
  shells---fundamentals}, Computational Fluid and Solid Mechanics, Springer,
  Heidelberg, second~ed., 2011.

\bibitem{DevlooFGG_13_ACC}
{\sc P.~R.~B. Devloo, A.~M. Farias, S.~M. Gomes, and J.~a.~L. Gon\c{c}alves},
  {\em Application of a combined continuous-discontinuous {G}alerkin finite
  element method for the solution of the {G}irkmann problem}, Comput. Math.
  Appl., 65 (2013), pp.~1786--1794.

\bibitem{Fluegge_60_SS}
{\sc W.~Fl\"{u}gge}, {\em Stresses in Shells}, Springer-Verlag,
  Berlin-G\"{o}ttingen-Heidelberg, 1960.
\newblock Translation of Statik und Dynamik der Schalen, Springer, 1934.

\bibitem{FrancisPT_90_BAA}
{\sc T.-L. Francis, V.~A. Pulmano, and D.~Thambiratnam}, {\em Bef analogy for
  axisymmetrically loaded cylindrical shells}, Computers \& Structures, 34
  (1990), pp.~281--285.

\bibitem{Girkmann_56_FEE}
{\sc K.~Girkmann}, {\em Fl\"{a}chentragwerke: {E}inf\"{u}hrung in die
  {E}lastostatik der {S}cheiben, {P}latten, {S}chalen und {F}altwerke},
  Springer-Verlag, Vienna, 1956.

\bibitem{Goldenveizer_61_TET}
{\sc A.~L. Gol'denve\u{\i}zer}, {\em Theory of elastic thin shells},
  International Series of Monographs in Aeronautics and Astronautics, Published
  for the American Society of Mechanical Engineers by Pergamon Press,
  Oxford-London-New York-Paris, 1961.
\newblock Translation from the Russian edited by G. Herrmann.

\bibitem{OlssonR_40_PBE}
{\sc R.~Gran~Olsson and E.~Reissner}, {\em A problem of buckling of elastic
  plates of variable thickness}, J. Math. Phys. Mass. Inst. Tech., 19 (1940),
  pp.~131--139.

\bibitem{HeuerK_17_RDM}
{\sc N.~Heuer and M.~Karkulik}, {\em A robust {DPG} method for singularly
  perturbed reaction-diffusion problems}, SIAM J. Numer. Anal., 55 (2017),
  pp.~1218--1242.

\bibitem{LinS_12_BFE}
{\sc R.~Lin and M.~Stynes}, {\em A balanced finite element method for
  singularly perturbed reaction-diffusion problems}, SIAM J. Numer. Anal., 50
  (2012), pp.~2729--2743.

\bibitem{Linss_10_LAM}
{\sc T.~Lin{\ss}}, {\em Layer-adapted meshes for reaction-convection-diffusion
  problems}, vol.~1985 of Lecture Notes in Mathematics, Springer-Verlag,
  Berlin, 2010.

\bibitem{MillerORS_96_FNM}
{\sc J.~J.~H. Miller, E.~O'Riordan, and G.~I. Shishkin}, {\em Fitted numerical
  methods for singular perturbation problems}, World Scientific Publishing Co.,
  Inc., River Edge, NJ, 1996.

\bibitem{Morley_76_ADS}
{\sc L.~S.~D. Morley}, {\em Analysis of developable shells with special
  reference to the finite element method and circular cylinders}, Philosophical
  Transactions of the Royal Society of London. Series A, Mathematical and
  Physical Sciences, 281 (1976), pp.~113--170.

\bibitem{Niemi_16_BCS}
{\sc A.~H. Niemi}, {\em Benchmark computations of stresses in a spherical dome
  with shell finite elements}, SIAM J. Sci. Comput., 38 (2016), pp.~B440--B457.

\bibitem{NiemiBPD_12_FEA}
{\sc A.~H. Niemi, I.~Babu\v{s}ka, J.~Pitk\"aranta, and L.~Demkowicz}, {\em
  Finite element analysis of the {Girkmann} problem using the modern hp-version
  and the classical h-version}, Engineering with Computers, 28 (2012),
  pp.~123--134.

\bibitem{OMalley_91_SPA}
{\sc R.~E. O'Malley, Jr.}, {\em Singular perturbations, asymptotic evaluation
  of integrals, and computational challenges}, in Asymptotic analysis and the
  numerical solution of partial differential equations ({A}rgonne, {IL}, 1990),
  vol.~130 of Lecture Notes in Pure and Appl. Math., Dekker, New York, 1991,
  pp.~3--16.

\bibitem{PitkaerantaBS_12_DRV}
{\sc J.~Pitk\"{a}ranta, I.~Babu\v{s}ka, and B.~Szab\'{o}}, {\em The dome and
  the ring: verification of an old mathematical model for the design of a
  stiffened shell roof}, Comput. Math. Appl., 64 (2012), pp.~48--72.

\bibitem{Reissner_37_RTB}
{\sc M.~E. Reissner}, {\em Remark on the theory of bending of plates of
  variable thickness}, Journal of Mathematics and Physics, 16 (1937),
  pp.~43--45.

\bibitem{RoosL_99_SCU}
{\sc H.-G. Roos and T.~Lin{\ss}}, {\em Sufficient conditions for uniform
  convergence on layer-adapted grids}, Computing, 63 (1999), pp.~27--45.

\bibitem{Schwab_98_php}
{\sc C.~Schwab}, {\em {$p$}- and {$hp$}-finite element methods}, Numerical
  Mathematics and Scientific Computation, The Clarendon Press, Oxford
  University Press, New York, 1998.

\bibitem{Shishkin_92_DAS}
{\sc G.~I. Shishkin}, {\em Discrete Approximation of Singularly Perturbed
  Elliptic and Parabolic Equations}, Russian Academy of Sciences, Ural Section,
  Ekaterinburg, 1992.
\newblock In Russian.

\bibitem{VentselK_01_TPS}
{\sc E.~Ventsel and T.~Krauthammer}, {\em Thin Plates and Shells}, CRC Press,
  New York, 2001.

\end{thebibliography}

\end{document}